\documentclass[1p]{article} 
\usepackage{amsmath,amsfonts,amsthm} 
\usepackage{subfigure}
\usepackage{latexsym}
\usepackage{amssymb}
\usepackage[dvips]{graphics}
\usepackage[dvips]{graphicx}
\newtheorem{theorem}{Theorem}

\newtheorem{remark}[theorem]{Remark}
\newtheorem{example}[theorem]{Example}

\newcommand{\R}{{\mathbb R}}
\newcommand{\Z}{{\mathbf Z}}
\newcommand{\Q}{{\mathbf Q}}

\renewcommand{\int}{\rm Int}
\newcommand{\E}{\mathbb E}

\newcommand{\PP}{\Bbb {P}}

\newcommand{\D}{{\mathfrak D}}

\newcommand{\p}{{\mathfrak p}}

\begin{document}

\title{Homological Domination in Large Random Simplicial Complexes}          
\author{A. Costa and M. Farber}        
\date{}          
\maketitle

\section{Introduction}

Manifolds and simplicial complexes traditionally appear in robotics as configuration spaces of mechanical systems. 
The classical approach becomes inadequate if we are dealing with {\it a large system}, i.e. with a system depending on a large number of metric parameters, since these parameters cannot be measured without errors and small errors make significant impact on the structure of the obtained space \cite{F}. A more realistic approach for modelling large systems is based on the assumption that their configuration spaces are random with properties described using the language of probability theory.

The mathematical study of topology of large random spaces started relatively recently and several different probabilistic models of random topological objects have appeared within the last 10 years, see \cite{CFK} and \cite{Ksurvey} for surveys. 
One may mention random surfaces \cite{PS}, random 3-dimensional manifolds
\cite{DT}, random configuration spaces of linkages
 \cite{F} and others.

The high dimensional large random simplicial complexes may also be used for modelling large networks, especially in situations 
when not only pairwise relations between the objects are important but also relations between triples, quadruples, etc.

 Linial, Meshulam and Wallach \cite{LM}, \cite{MW} initiated
an important analogue of the classical Erd\H os--R\'enyi model \cite{ER} of random graphs in the situation of high-dimensional simplicial complexes. 
The random simplicial complexes of \cite{LM}, \cite{MW} are $d$-dimensional, have the complete $(d-1)$-skeleton and their randomness shows
only in the top dimension. Some interesting results about the topology of random 2-complexes in the Linial--Meshulam model were obtained in \cite{BHK}, \cite{CCFK}, \cite{CF1}. 


A different model of random simplicial complexes was studied by M. Kahle \cite{Kahle1} and by some other authors, see for example \cite{CFH}.
These are clique complexes of random Erd\H os--R\'enyi graphs, i.e. here one starts with a random graph in the 
Erd\H os--R\'enyi model and declares as a simplex every subset of vertices which form a {\it clique} (a subset such that every two vertices are connected by an edge). 
Compared with the Linial - Meshulam model, the clique complex has randomness in dimension one but it influences its structure in all the higher dimensions.

In \cite{CF14}, \cite{CF15} the authors initiated the study of a more general and more flexible model of random simplicial complexes with randomness in all dimensions. 
Here one starts with a set of $n$ vertices and retain each of them with probability $p_0$; on the next step one connects every pair of retained vertices by an edge with probability $p_1$, and then fills in every triangle in the obtained random graph with probability $p_2$, and so on. 
As the result we obtain a random simplicial complex depending on the set of probability parameters 
$$\p = (p_0, p_1, \dots, p_r), \quad 0\le p_i\le 1.$$
Our multi-parameter random simplicial complex includes both Linial-Meshulam and random clique complexes as special cases. 
The topological and geometric properties of multi-parameter random simplicial complexes depend on the whole set of parameters and their thresholds can be understood as boundaries of convex subsets and not as single numbers as in all the previously studied models. 

In this paper we state {\it the homological domination principle} (Theorem \ref{thmdom}) for random simplicial complexes, claiming that the Betti number 
in one specific dimension $k=k(\p)$ (which is explicitly determined by the probability multi-parameter $\p$) significantly dominates the Betti numbers in all other dimensions. 
We also state and discuss evidence for two interesting conjectures which strengthen the homological domination principle; they claim that homology in dimensions below $k=k(\p)$ vanishes and homology in dimensions above $k=k(\p)$ is generated by bounding degree zero faces. 
 \vskip 0.2cm 
 This research was supported by the EPSRC research council. 




\section{Random simplicial complexes depending on several probability parameters.}  

\subsection{The model} Let $\Delta_n$ denote the simplex with the vertex set $\{1, 2, \dots, n\}$. 
We view $\Delta_n$ as an abstract simplicial complex of dimension $n-1$. 
For a simplicial subcomplex $Y\subset \Delta_n$, we denote by $f_i(Y)$ the number of {\it $i$-faces} of $Y$ (i.e. $i$-dimensional simplexes of $\Delta_n$ contained in $Y$). 
An external face of a subcomplex $Y\subset \Delta_n$ is a simplex $\sigma \subset \Delta_n$ such that $\sigma \not\subset Y$ but the boundary of $\sigma$ is contained in $Y$, 
$\partial \sigma \subset Y$. The symbol $e_i(Y)$ will denote the number of 
$i$-dimensional external faces of $Y$. 



%

Fix an integer $r\ge 0$ and a sequence ${\mathfrak p}=(p_0, p_1, \dots, p_r)$ of real numbers satisfying $0\le  p_i\le 1.$ Denote  
$q_i=1-p_i.$
We consider the probability space ${\Omega_n^r}$ consisting of all subcomplexes 
$Y\subset \Delta_n,$ with $ \dim Y\le r.$
The probability function
\begin{eqnarray}
\PP_{r,\p}: {\Omega_n^r}\to \R
\end{eqnarray} is given by the formula
\begin{eqnarray}\label{def1}
\PP_{r, \p}(Y) \, 
&=& \, \prod_{i=0}^r p_i^{f_i(Y)}\cdot \prod_{i=0}^r q_i^{e_i(Y)}
\end{eqnarray}
In (\ref{def1}) we use the convention $0^0=1$; in other words, if $p_i=0$ and $f_i(Y)=0$ then the corresponding factor in (\ref{def1}) equals 1; similarly if some $q_i=0$ and $e_i(Y)=0$. 
One may show that $\PP_{r, \p}$ is indeed a probability function, i.e. 
$
\sum_{Y\subset \Delta_n^{(r)}}\PP_{r, \p}(Y) = 1, 
$
see \cite{CF14}.



\subsection{Important special cases.} The multi-parameter model we consider in this paper turns into some important well known models in several special cases:

When $r=1$ and $\p=(1, p)$ we obtain the classical model of random graphs of Erd\H os and R\'enyi \cite{ER}. 

When $r=2$ and $\p= (1,1,p)$ we obtain the Linial - Meshulam model of random 2-complexes \cite{LM}. 

When $r$ is arbitrary and fixed and $\p = (1, 1, \dots, 1, p)$ we obtain the random simplicial complexes of Meshulam and Wallach \cite{MW}. 

For $r=n-1$ and $\p=(1,p, 1,1, \dots, 1)$ one obtains the clique complexes of random graphs studied in \cite{Kahle1}.

\subsection{Connectivity and simple-connectivity}

We state below (in a simplified form) a result obtained in \cite{CF15}: 

\begin{theorem} \label{thm1} Consider a random simplicial complex $Y\in \Omega_n^r$ with respect to the probability measure $\PP_{r, \p}$ where 
$\p=(p_0, \dots, p_r)$ and $p_i=n^{-\alpha_i}$ for $i=0, \dots, r$. We assume that the numbers $\alpha_i\ge 0$ do not depend on $n$. Then for 
$$r\ge 1, \quad \quad \alpha_0+\alpha_1<1$$
a random complex $Y$ is connected, a.a.s. Besides, for 
$$r\ge 2, \quad \quad\alpha_0+3\alpha_1 +2\alpha_2 <1$$
a random complex $Y$ is simply-connected, a.a.s. 
\end{theorem} 

See \cite{CF15} for proofs and more details. 


\section{Homological domination principle.}

Consider the following linear functions:
\begin{align*}
\psi_k(\alpha)=\sum_{i=0}^r \binom{k}{i}\alpha_i,\quad\quad k=0,\ldots, r, 
\end{align*}
where $$\alpha=(\alpha_0, \dots, \alpha_r)\in \R^{r+1}.$$ 
We use the conventions that $\binom{k}{i}=0$ for $i>k$ and $\binom{0}{0}=1$.
We shall assume that $\alpha\in \R_+^{r+1}$, i.e. $\alpha=(\alpha_0, \dots, \alpha_r)$ with $\alpha_i\ge 0$ for all $i$. 
Since $\binom k i < \binom {k+1} i$ for $i>0$ we see that 
\begin{align}\label{ineq1}
\psi_0(\alpha)\leq\psi_1(\alpha)\leq\psi_2(\alpha)\leq\ldots\leq\psi_r(\alpha).
\end{align}
Moreover, if for some  $j\geq 0$ one has $\psi_j(\alpha)<\psi_{j+1}(\alpha)$ then
\begin{align}\label{ineq2}\psi_j(\alpha)<\psi_{j+1}(\alpha)<\psi_{j+2}(\alpha)<\ldots<\psi_r(\alpha).
\end{align}

Define the following convex domains (open sets) in $\R^{r+1}_+$:
\begin{eqnarray}
{\mathfrak D}_k = \{\alpha\in \R^{r+1}_+\, ;\,  \psi_k(\alpha)<1<\psi_{k+1}(\alpha)\}, \quad k=0, 1, \dots, r-1.
\end{eqnarray}
 \begin{figure}[h]
\centering
\includegraphics[width=0.65\textwidth]{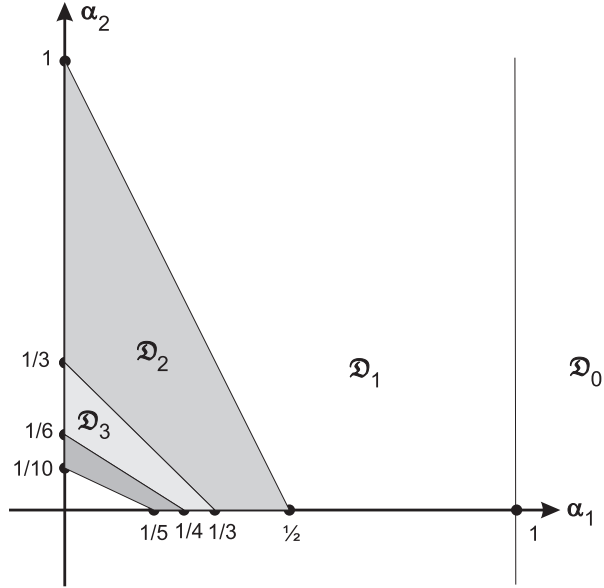}
\caption{Intersections of the domains $\D_k$ with the plane $\alpha_0=\alpha_3=\dots=0$. }\label{domains}
\end{figure}
One may also introduce the domains 
$$\D_{-1}=  \{\alpha\in \R^{r+1}_+\, ;\,  1<\psi_0(\alpha)\}, \quad \D_{r}=  \{\alpha\in \R^{r+1}_+\, ;\,  \psi_r(\alpha)<1\}.$$
The domains $\D_{-1}, \D_0, \D_1,\dots, \D_r$ are disjoint and their union is $$\bigcup_{j=-1}^r \D_j \, = \, \R^{r+1}_+- \bigcup_{i=0}^r H_i$$
where $H_i$ denotes the hyperplane given by the linear equation $\psi_i(\alpha)=1$. We shall see that each hyperplane $H_i$ 
correspond to {\it phase transitions in homology}; if $\alpha\in \D_{-1}$ then the random complex $Y$ is $\emptyset$, a.a.s.
Conjecturally (see \S \ref{conj} below), when $\alpha$ crosses the hyperplane $H_j$ and moves from the domain $\D_{j-1}$ to the domain $\D_{j}$, the random complex $Y$ changes from being homotopically $(j-1)$-dimensional and becoming homotopically $j$-dimensional.

Next, we define a non-negative quantity 
\begin{eqnarray}\label{ealpha}
e(\alpha) = \min_k \{|1-\psi_k(\alpha)|\}. 
\end{eqnarray}
Note that $e(\alpha) = \min\{1-\psi_k(\alpha), \psi_{k+1}(\alpha) -1\}>0$
assuming that $\alpha\in {\mathfrak D}_k$. 

\begin{theorem}\label{thmdom} Consider a multi-parameter random simplicial complex $Y\in \Omega_n^r$ with respect to the probability measure 
$$\PP_{r, \p}:\Omega_n^r\to \R,\quad \mbox{where } \quad 
\p = n^{-\alpha},$$ i.e. $$\p=(p_0, p_1, \dots, p_r), \qquad p_i=n^{-\alpha_i}.$$
The Betti numbers $b_j: \Omega_n^r\to \R$ of various dimensions $j=0, 1, \dots, r$ are random variables and for $\alpha \in {\mathfrak D}_k$, 
where $k=0, \dots , r$, their expectations have the following properties:
\begin{enumerate}
\item
Firstly, for $\alpha\in {\mathfrak D}_k$ and $n$ large enough one has
\begin{eqnarray}
\E(b_k) \ge \frac{n^{e(\alpha)}}{(r+1)!}\cdot \E(b_j), \qquad j\not= k. 
\end{eqnarray}
In other words, for $\alpha \in {\mathfrak D}_k$ the $k$-th Betti number significantly dominates all other Betti numbers. 
\item Secondly, 
for $\alpha \in {\mathfrak D}_k$, 
\begin{eqnarray}
\E(b_k) \sim \frac{1}{(k+1)!} \cdot n^{k+1-\sum_{i=0}^k \psi_i(\alpha)}.
\end{eqnarray}
\end{enumerate}
\end{theorem}

A detailed proof of Theorem \ref{thmdom} will be published elsewhere; it is based on a Morse theory arguments and on the following fact established in \cite{CF14}. Fix $0\le \ell \le r$ and consider the number of simplexes of dimension $\ell$ as a random variable 
$f_\ell:\Omega_n^r\to \R.$
Then the expectation of $f_\ell$ is 
$$\E(f_\ell) \, = \,  \binom n {\ell+1} \cdot \prod_{i=0}^\ell p_i^{\binom {\ell+1} {i+1}}$$
and for $p_i=n^{-\alpha_i}$ one has 
$$\E(f_\ell) \sim \frac{1}{(\ell+1)!}\cdot n^{\ell+1-\sum_{i=0}^\ell \psi_i(\alpha)}.$$
There is also a stronger version of Theorem \ref{thmdom} which operates with the Betti numbers $b_j$ rather than with their expectations.

\section{Dimension}

In this section we clarify the role of the domains $\D_k\subset \R^{r+1}_+$ play regarding the issue of the dimension of a random simplicial complex. 

Consider the following map 
$\Phi: \R^{r+1}_+\to \R^{r+1}_+$
given by 
$$\Phi(\alpha_0, \alpha_1, \dots, \alpha_r) = (\alpha'_0, \alpha'_1, \dots, \alpha'_r), $$
where $$\quad\quad  \alpha'_i \, =\,   \frac{\alpha_i}{i+1}, \quad i=0, \dots, r.$$
\begin{theorem}\label{thmdim} 
Consider the probability measure $\PP_{r, \p}$ on $\Omega_n^r$ where 
$\p = n^{-\alpha}$ with $\alpha\in \R^{r+1}_+$. If $\Phi(\alpha)\in \R^{r+1}_+$ lies in $\D_k$, i.e. $\Phi(\alpha)\in \D_k$, then the dimension of a random complex $Y\in \Omega_n^r$ equals $k$, a.a.s.
\end{theorem}

Hence we see that the pre-images of the domains $\D_k$ under the map $\Phi$ are exactly the domains where a random complex has geometric dimension $k$. 

A proof of Theorem \ref{thmdim} can be found in \cite{CF14}.

\section{The General Picture}\label{conj}

Consider a random simplicial complex $Y\in \Omega_n^r$ with respect to the probability measure $\PP_{r, \p}:\Omega_n^r\to \R,$ where 
$\p = n^{-\alpha},$ i.e. $\p=(p_0, p_1, \dots, p_r)$, $p_i=n^{-\alpha_i},$ and $\alpha=(\alpha_0, \dots, \alpha_r)\in \R^{r+1}_+$. 
In this section we state two interesting conjectural statements related to the properties of multi-parameter random simplicial complexes. 
We examine these statements in several special cases and show that they are consistent with known results. 

A simplex $\sigma\subset Y$ of a simplicial complex has {\it degree zero} if it is not a face of any simplex of higher dimension. 
A degree zero simplex $\sigma\subset Y$ is said to be {\it bounding} if $\partial \sigma$ bounds a chain with $\Q$ coefficients in $Y-\int(\sigma)$. 
Removing a bounding degree zero simplex of dimension $d$ reduces the $d$-dimensional Betti number by $1$ and does not affect the other Betti numbers. 

\vskip 0.2cm

We believe that the following statements are true:
\vskip 0.2cm

{\bf  A}: 
 {\it For $\alpha\in \D_k$, consider a random complex $Y\in \Omega_n^r$, and 
 let $Y'\subset Y$ be obtained from $Y$ by removing a set of degree zero bounding simplexes of dimension $>k$. Then 
 $Y'$ collapses simplicity to a $k$-dimensional subcomplex 
$Z\subset Y$, a.a.s. In particular, the homology groups of $Y'$ with any coefficients vanish in all dimensions $>k$, a.a.s. 
}

\vskip 0.2cm
{\bf B}: 
{\it For $\alpha\in \D_k$, the reduced homology of the random complex $Y\in \Omega_n^r$ with integral coefficients vanishes in all dimensions less than $k$, a.a.s.
}
\vskip 0.2cm
A weaker version of statement {\bf B} dealing with homology with {\it rational} coefficients can be analysed using the well known Garland's method (see for example, \cite{Kahle3}) based on spectral analysis of the combinatorial Laplacian of links of simplexes. 
%
%
%
%
%
%
%
%
%
\vskip 0.2cm 
Next we consider special cases:
\vskip 0.2cm 
{\bf Case I}. 
Consider the special case when $r=2$ and $\alpha=(0,0, \alpha_2)$; this case corresponds to the Linial - Meshulam model of random 2-complexes \cite{LM}. 
In this case the domains $\D_{-1}=\D_0=\emptyset$ are empty and the other domains are 
$$\D_1=\{\alpha_2>1\}, \quad \mbox{and}\quad  \D_2=\{\alpha_2<1\}.$$
Statement {\bf A} in this special case follows from Theorem 1 from \cite{CCFK} which states that for $\alpha\in \D_1$ the random complex collapses to a graph; if $\alpha\in \D_2$ there is nothing to prove. 
Statement {\bf B} in this special case reduces to the statement that for $\alpha\in \D_2$ the reduced integral homology of a random complex vanishes 
in dimension $1$; this is exactly the content of the main theorem from \cite{KHP}. 

\vskip 0.2cm 
{\bf Case II}. Consider now the case when $\alpha=(0, \alpha_1, 0, 0\dots, 0)$ which corresponds to the $r$-skeleton of the clique complex of a random 
Erd\H os -- R\'enyi graph with edge probability $p=n^{-\alpha_1}$. The domains $\D_i$ in this case are $\D_{-1}=\emptyset$, 
$\D_0= \{\alpha_1>1\}, $
$$\D_i=\left\{\frac{1}{i+1}< \alpha_1< \frac{1}{i}\right\}, \quad \mbox{for}\quad i=1, \dots, r-1$$
and 
$\D_r =\{0\le \alpha_1<1/r\}$. It is known that for $\alpha\in\D_i$, firstly, 
the integral homology of the random complex $H_k(Y;\Z)=0$ vanish in dimensions $k>i$
(see \cite{Kahle1}, Theorems 3.6) and, secondly, the reduced homology groups $\tilde H_j(Y;\Q)$ with rational coefficients vanish in all dimensions $j<i$, see \cite{Kahle3}, Theorem 1.1. Both statements would follow from {\bf A} and {\bf B}. Besides, Theorem A from \cite{CFH} implies that for $\alpha\in \D_1$ the random complex collapses simplicity to a graph, which is consistent with Statement {\bf A}.

\bibliographystyle{amsalpha}

\end{document}